\documentclass[nobysame,initials]{amsart}
\usepackage{xcolor}
\usepackage[normalem]{ulem}
\usepackage{bookmark}
\usepackage[alphabetic]{amsrefs}
\usepackage[makeroom]{cancel}
\usepackage{bbold}
\usepackage{enumerate}
\usepackage{fullpage}

\allowdisplaybreaks
\newcommand{\R}{\mathbb{R}}
\newcommand{\Red}{\R^d}
\newcommand{\bd}{\mathrm{bd}}
\newcommand{\E}{\mathbb{E}}
\newcommand{\V}{\mathrm{Vol}}
\newcommand{\PP}{\mathbb{P}}

\newcommand{\id}{\mathrm{Id}}
\newcommand{\dd}{\;\mathrm{d}}
\newcommand{\supp}{\mathrm{supp}}
\newcommand{\iid}{i=1,\ldots, d}
\newcommand{\di}{\;\mathrm d}
\newcommand{\st}{:\;}
\newcommand{\spann}{\mathrm{span}}
\newcommand{\rank}{\mathrm{rank}}
\newcommand{\mm}{\overline{m}}
\newcommand{\dbound}{d(d+1)/2}
\newtheorem{theorem}{Theorem}
\newtheorem{proposition}{Proposition}
\newtheorem{lemma}{Lemma}

\theoremstyle{definition}

\title{On the volume bound in the Dvoretzky--Rogers lemma}

\author{Ferenc Fodor}
\address{Department of Geometry, Bolyai Institute, University of Szeged, Aradi 
v\'ertan\'uk tere 1, 6720 Szeged, Hungary}
\email{fodorf@math.u-szeged.hu}

\author{M\'arton Nasz\'odi}
\address{Department of Geometry, E\"otv\"os University, Budapest, Hungary}
\email{marton.naszodi@math.elte.hu}

\author{Tam\'as Zarn\'ocz}
\address{Bolyai Institute, University of Szeged, Aradi 
v\'ertan\'uk tere 1, 6720 Szeged, Hungary}
\email{tzarnocz@math.u-szeged.hu}



\begin{document}

\begin{abstract}
The classical Dvoretzky--Rogers lemma provides a deterministic algorithm by 
which, from any set of isotropic vectors in Euclidean $d$-space, one can select 
a subset of $d$ vectors whose determinant is not too small. 
Subsequently, Pelczy\'nski and Szarek improved this lower bound by a 
factor depending on the dimension and the number of vectors.

Pivovarov, on the other hand, 
determined the expectation of the square of the volume of parallelotopes 
spanned by $d$ independent random vectors in $\R^d$, each one chosen according 
to an isotropic measure.
We extend Pivovarov's result to a class of more general probability 
measures, which yields that the volume bound in the 
Dvoretzky--Rogers lemma is, in fact, equal to the expectation of the squared 
volume of random parallelotopes spanned by isotropic vectors. This allows us to 
give a probabilistic proof of the improvement of Pelczy\'nski and Szarek, and 
provide a lower bound for the probability that the volume of such a random 
parallelotope is large.
\end{abstract}

\maketitle

\section{Introduction}

Given a set of isotropic vectors in Euclidean $d$-space $\R^d$ (see definition 
below), the \emph{Dvoretzky--Rogers lemma} states that one may select a subset 
of $d$ ``well spread out'' vectors. As a consequence, the determinant of these 
$d$ vectors is at least $\sqrt{d!/d^d}$. This selection is deterministic: we 
start with an arbitrary element of the set, and then select more vectors 
one-by-one in a certain greedy manner.

Pivovarov \cite{P2010}*{Lemma~3, p. 49}, on the other hand, chooses $d$ 
vectors randomly and then computes the expectation of the square 
of the resulting determinant. In this note, we 
extend Pivovarov's result to a wider class of measures, and 
apply this extension to obtain the improved lower 
bound of Pelczy\'nski and Szarek, cf. \cite{PSz}{Proposition~2.1}, on the 
maximum of the volume of parallelotopes spanned by 
$d$ vectors from the support of the measure. Thus, we give a 
probabilistic interpretation of the volume bound 
in the Dvoretzky--Rogers lemma.

We denote the Euclidean scalar product by $\langle\cdot, \cdot\rangle$, the 
induced norm by $|\cdot |$. We use the usual notation $B^d$ for the unit ball 
of $\R^d$ centered at the origin $o$, and $S^{d-1}$ for its boundary $\bd B^d$. 
We call a compact convex set 
$K\subset \R^d$ with non-empty interior a \emph{convex body}. For detailed 
information on the properties of convex bodies, we refer to the books by Gruber 
\cite{Gru07} and Schneider \cite{Sch14}. 

Let $\id_d$ be the identity map on $\Red$.
For $u,v\in \R^d$, let $u\otimes v:\Red\rightarrow\Red$ denote the \emph{tensor 
product} of $u$ and $v$, that is, $(u\otimes v)(x)=\langle v, x\rangle u$ for 
any $x\in\Red$. Note that when $u\in S^{d-1}$ is a unit vector, $u\otimes u$ is 
the orthogonal projection to the linear subspace spanned by $u$. 

For two functions $f(n), g(n)$, we use the notation $f(n)\sim 
g(n)$ (as $n\to\infty$) if $\lim_{n\to\infty}f(n)/g(n)=1$.

An \emph{isotropic measure} is a probability measure $\mu$ on $\R^d$ with the 
following two properties.

\begin{equation}\label{eq:isotropictesnor}
\int_{\R^d} x\otimes x \dd \mu(x)=\id_d,
\end{equation}
and the center of mass of $\mu$ is at the origin, that is,
\begin{equation}\label{eq:isotropiccenterofmass}
\int_{\R^d} x \dd \mu(x)=0.
\end{equation}

Pivovarov~\cite{P2010} proved the following statement about the 
volume of random parallelotopes spanned by $d$ independent, isotropic vectors.

\begin{lemma}[Pivovarov~\cite{P2010}, Lemma~3]\label{lem:pivovarov}
	Let $x_1,\ldots, x_d$ be independent random vectors distributed 
	according to 
	the isotropic measures $\mu_1, \ldots, \mu_d$ in $\R^d$.
	Assume that $x_1,\ldots, x_d$ are linearly independent with probability 
$1$. 
	Then 
	\begin{equation}\label{eq:Pivovarov}
	\E([\det (x_1,\ldots, x_d)]^2)=d!.
	\end{equation}
\end{lemma}

We note that Lutwak, Yang and Zhang in \cite{LYZ04}*{\S 2} established 
similar results for the case of discrete isotropic measures, which could also 
be used to prove the volumetric bounds in Theorem~\ref{main}, see, for 
example, \cite{LYZ04}*{formula (2.5) on page 167}.

We extend Lemma~\ref{lem:pivovarov} to a more general class of measures
in the following way.

\begin{lemma}\label{mainmainlemma}
	Let $x_1,\ldots, x_d$ be independent random vectors distributed 
according to the probability measures $\mu_1, \ldots, \mu_d$ in $\R^d$ satisfying 
	\eqref{eq:isotropictesnor}. Assume that $\mu_i(\{0\})=0$ for \iid. Then 
	\eqref{eq:Pivovarov} holds.
\end{lemma}

We provide a simple and direct proof of Lemma~\ref{mainmainlemma} in 
Section~\ref{sec:mainlemma}.

Lemmas~\ref{lem:pivovarov} and \ref{mainmainlemma} yield the value 
of the second moment of the volume of random parallelotopes with isotropic 
generating vectors. 
On the other hand, Milman and Pajor  \cite{MP89}*{\S3.7} gave a lower bound 
for the $p$-th moment (with $0<p<2$) of this volume in the case when the 
generating vectors are selected according to the uniform distribution 
from an isotropic and origin-symmetric convex body; for more general results, 
cf. 
\cite{BGVV}*{\S3.5.1}. All of the previously mentioned results hold in 
\emph{expectation}.

As a different approach, we mention Pivovarov's work \cite{P2010}, where lower 
bounds on the volume of a random parallelotope are shown to hold \emph{with 
high probability} under the assumption that the measures are log-concave.

For more information on properties of random parallelotopes, and random 
polytopes in general, we refer to the book by Schneider and Weil \cite{SW08}, 
the survey by Schneider \cite{Sch17}, and the references therein.

In this paper, our primary, geometric motivation in studying isotropic 
measures is the 
following celebrated theorem of John \cite{F48}, which we state in the refined 
form obtained by Ball \cite{B92} (see also \cite{B97}).

\begin{theorem}\label{John}
	Let $K$ be a convex body in $\R^d$. Then there exists a unique 
ellipsoid of 
	maximal volume contained in $K$. Moreover, this maximal volume 
ellipsoid is 
the 
	$d$-dimensional unit ball $B^d$ if and only if there exist vectors 
$u_1,\ldots, 
	u_m\in \bd K\cap S^{d-1}$ and (positive) real numbers $c_1,\ldots, 
c_m>0$ 
such 
	that 
	\begin{equation}\label{isotropy}
	\sum_{i=1}^{m}c_i u_i\otimes u_i=\mathrm{Id}_d,
	\end{equation}
	and
	\begin{equation}\label{centeredness}
	\sum_{i=1}^{m}c_iu_i=0.
	\end{equation} 
\end{theorem}

Note that taking the trace in \eqref{isotropy} yields $\sum_{i=1}^{m}c_i=d$. 
Thus, 
the Borel measure $\mu_K$ on $\sqrt{d}S^{d-1}$ with 
$\supp \mu_K=\{\sqrt{d}u_1,\ldots, \sqrt{d}u_m\}$ and 
$\mu_K(\{\sqrt{d}u_i\})=c_i/d$ $(i=1,\ldots, m)$ 
is a discrete isotropic measure.  

If a finite system of unit vectors $u_1,\ldots, u_m $ in $\R^d$, together with 
a set 
of positive weights $c_1,\ldots, c_m$ satisfies  \eqref{isotropy} and 
\eqref{centeredness}, then we say that it forms a \emph{John decomposition of 
	the identity}. For each convex body $K$, there exists an affine image 
$K^{\prime}$ of $K$ for which the maximal volume ellipsoid contained in 
$K^{\prime}$ is $B^d$, and $K^{\prime}$ is unique up to orthogonal 
transformations of $\Red$.

The classical lemma of Dvoretzky and Rogers \cite{DR1950} states that in a John 
decomposition of the identity, one can 
always find $d$ vectors such that the selected vectors are not too far from an 
orthonormal system.

\begin{lemma}[Dvoretzky--Rogers lemma \cite{DR1950}]
	Let $u_1,\ldots, u_m\in S^{d-1}$ and $c_1,\ldots, c_m>0$ such that 
	\eqref{isotropy} holds. Then there exists an orthonormal basis 
$b_1,\ldots, 
	b_d$ of $\Red$ and a subset $\{ x_1,\ldots, 
	x_d\}\subset \{u_1,\ldots, u_m\}$ with $x_j\in\mathrm{lin}\{b_1,\ldots, 
	b_j\}$ and
	\begin{equation}
	\sqrt{\frac{d-j-1}{d}}\leq \langle x_j, b_j\rangle\leq 1
	\end{equation}
	for $j=1,\ldots, d$. 
\end{lemma}

Consider the parallelotope $P$ spanned by the selected $d$ vectors 
$x_1,\ldots, x_d$. The volume of $P$ is bounded from below by
\begin{equation}\label{eq:volboundDvoRo}
(\V(P))^2= [\det (x_1,\ldots, x_d)]^2\geq \frac{d!}{d^d}.
\end{equation}

Our study of \eqref{eq:volboundDvoRo} is motivated in part by the recent proof 
\cite{Nas} of a conjecture of B\'ar\'any, Katchalski and Pach, where this bound 
is heavily relied on.

The main results of this paper are the following two statements. 
Theorem~\ref{main} is essentialy the same as Proposition~2.1 of Pelczy\'nski 
and 
Szarek \cite{PSz}, however, here we give a probabilistic proof and 
interpretation. In Theorem~\ref{main} (\ref{csvaltozat}) and (\ref{main3}), we 
also note that when $m$ is 
small the improvement on the original Dvoretzky--Rogers bound is 
larger.

\begin{theorem}\label{main}
	Let $u_1,\ldots, u_m\in S^{d-1}$ be unit vectors 
	satisfying \eqref{isotropy} with some $c_1,\ldots,\linebreak[0]c_m > 0$.
	Then there is a subset $\{ x_1,\ldots, x_d\}\subset \{u_1,\ldots, 
u_m\}$ 
with
		\begin{equation*}
		[\det (x_1,\ldots, x_d)]^2 \geq \gamma(d,\mm)\cdot 
\frac{d!}{d^d},
		\end{equation*}
where $\gamma(d,\mm)=\frac{\mm^d}{d!}\binom{\mm}{d}^{-1}$, 
and $\mm=\min\{m,\dbound\}$. 
		
Moreover, for $\gamma(d,\mm)$, we have
	\begin{enumerate}[(i)]
	  \item \label{main1}
		$\gamma(d,\mm)\geq \gamma(d,\dbound)\geq 3/2$ for any $d\geq 2$ 
and $m\geq 
d$. And $\gamma(d,\dbound)$ is monotonically increasing, and 
$\lim_{d\to\infty}\gamma(d,\dbound)=e$.
	  \item \label{csvaltozat} Fix a $c>1$, and consider the case when 
$m\leq cd$ 
with $c\geq1+1/d$. Then
		\begin{equation*}
		\gamma(d,m)\geq \gamma(d,\lceil cd\rceil)\sim 
\sqrt{\frac{c-1}c}\left( 
	  \frac{c-1}c\right)^{(c-1)d}e^d,\quad \text{as } d\to\infty.
		\end{equation*}
	  \item \label{main3} Fix an integer $k\geq1$, and consider the case 
when 
$m\leq d+k$. Then
	  	\begin{equation*}
	  \gamma(d,m)\geq\gamma(d,d+k)\sim 
\frac{k!e^k}{\sqrt{2\pi}}\frac{e^d}{(d+k)^{k+1/2}}, \text{ as } d\to\infty.
	  	\end{equation*}
	\end{enumerate}
\end{theorem}

We note that in (\ref{csvaltozat}) and (\ref{main3}), the 
improvements are exponentially large in $d$ as $d$ tends to infinity. 

The following statement provides a lower bound on the probability 
that $d$ independent, identically distributed 
random vectors selected from $\{u_1,\ldots, u_m\}$ according to 
the distribution determined by the weights $\{c_1,\ldots, c_m\}$ has large 
volume.

\begin{proposition}\label{mainProbability}
Let $\lambda\in(0,1)$. With the notations and assumptions of 
Theorem~\ref{main}, if we choose the vectors $x_1,\ldots,x_d$
independently according to the distribution $\PP(x_{\ell}=u_i)=c_i/d$ for each 
$\ell=1,\ldots,d$ and $i=1,\ldots,m$, then 
with probability at least $(1-\lambda)e^{-d}$, we have that 
\[[\det (x_1,\ldots, x_d)]^2 \geq \lambda\gamma(d,\mm)\cdot \frac{d!}{d^d}.\]
\end{proposition}

The geometric interpretation of Theorem~\ref{main} is the following. If $K$ is 
a convex polytope with $n$ facets, and $B^d$ is the maximal volume ellipsoid in 
$K$, then the number of contact points $u_1,\ldots,u_m$ in John's theorem is at 
most $m\leq n$. Thus, Theorem~\ref{main} yields a simplex in $K$ of not too 
small volume, with one vertex at the origin.

In particular, consider $k=1$ in Theorem~\ref{main} \eqref{main3}, that is, 
when $K$ is the regular simplex whose inscribed ball is $B^d$. Then the John 
decomposition of the identity determined by $K$ consists of $d+1$ unit 
vectors that determine the vertices of a regular 
$d$-simplex inscribed in $B^d$, which we denote by $\Delta_d$, and note that
$
\V(\Delta_d)=(d+1)^{\frac{d+1}{2}} / (d^{d/2}d!).
$
Clearly, in this John decomposition of the identity, the volume of the simplex 
determined by any $d$ of the vectors $u_1,\ldots, u_{d+1}$ is 
\begin{equation}\label{eq:boundforsimplex}
\V(\Delta_d)/(d+1)=\frac{(d+1)^{\frac{d-1}{2}}}{d^{d/2}d!}.
\end{equation}
By Theorem~\ref{main}, we obtain that
\[\max [\det (u_{i_1},\ldots, u_{i_d})]^2{\geq}
\frac{(d+1)^{d-1}}{d!}\cdot \frac{d!}{d^d}=\frac{(d+1)^{d-1}}{d^d},
\]
which yields the same bound for the largest volume simplex as the 
right-hand-side of \eqref{eq:boundforsimplex}. Thus, 
Theorem~\ref{main} is sharp in this case.


We will use the following theorem in our argument.
\begin{theorem}[\cites{F48, P, B92, GS2005}]\label{thm:GS}
If a set of unit vectors satisfies \eqref{isotropy} (resp., \eqref{isotropy} 
and \eqref{centeredness}) with some positive scalars $c'_i$, 
then a subset of 
$m$ elements also satisfies \eqref{isotropy} (resp., \eqref{isotropy} and 
\eqref{centeredness}) with some positive scalars $c_i$, where 
\begin{equation}\label{eq:GruberSchuster} d+1\leq m\leq d(d+1)/2 \end{equation} 
(resp., $d+1\leq m\leq d(d+3)/2$).
\end{theorem}

In Section~\ref{sec:GSproof}, we outline a proof of Theorem~\ref{thm:GS} for 
two reasons. First, we will use the part when only \eqref{isotropy} is 
assumed, which is only implicitly present in \cite{GS2005}. 
Second, 
in \cite{GS2005}, the result is described in terms of the contact points of a 
convex body with its maximal volume ellipsoid, that is, in the context of 
John's theorem. We, on the other hand, would like to give a presentation where 
the linear algebraic fact and its use in convex geometry are separated. 
Nevertheless, our proof is very close to the one given in \cite{GS2005}.

\section{Proof of Lemma~\ref{mainmainlemma}}\label{sec:mainlemma}
The idea of the proof is to slightly rotate each distribution so that the 
probability that the $d$ vectors are linearly independent is 
$1$. Then we may 
apply Pivovarov's lemma, and use a limit argument as the $d$ 
rotations each tend to the identity.

Let $A_1,\ldots,A_d$ be matrices in $SO(d)$ chosen independently of each other 
and of the $x_i$s according to the unique Haar probability measure on $SO(d)$. 
Fix an arbitrary non-zero 
unit vector $e$ in $\R^d$. Note that $A_ix_i/|x_i|$ and $A_ie$ have the same 
distribution: both are uniformly chosen points of the unit sphere according to 
the uniform probability distribution on $S^{d-1}$. A bit more 
is true: the joint distribution of $A_1x_1/|x_1|,\ldots,A_dx_d/|x_d|$ and the 
joint distribution of $A_1e,\ldots,A_de$ are the same: they are independently 
chosen, uniformly distributed points on the unit sphere.
It follows that

\begin{equation*}
\PP(A_1x_1,\ldots,A_dx_d \text{ are lin. indep.})=
\PP(A_1e,\ldots,A_de \text{ are lin. indep.})=1.
\end{equation*}
Denote the Haar measure on $Z:=SO(d)^d$ by $\nu$. Thus, we have
\begin{equation*}
1=\PP(A_1x_1,\ldots,A_dx_d \text{ are lin. indep.})=
\end{equation*}
\begin{eqnarray*}
\int_Z\int_{\R^d}\int_{\R^d}\ldots\int_{\R^d}\mathbb1_{\{A_1x_1,\ldots,A_dx_d 
		\text{ are lin. indep.}\}}(x_1,\ldots,x_d,A_1,\ldots,A_d)
	\\
	\di \mu_1(x_1)\ldots\di\mu_d(x_d)\di \nu(A_1,\ldots,A_d)
\end{eqnarray*}
\begin{equation*}
= \int_Z \PP(A_1x_1,\ldots,A_dx_d \text{ are lin. 
	indep.}\, |\, A_1,\ldots,A_d) \di \nu(A_1,\ldots,A_d),
\end{equation*}
where $\mathbb1$ denotes the indicator function.

Thus,
\begin{equation}\label{eq:allchoicesgood}
1=\PP\bigg[\PP(A_1x_1,\ldots,A_dx_d \text{ are lin. 
	indep.}\, |\, A_1,\ldots,A_d)=1\bigg].
\end{equation}
We call a $d$-tuple $(A_1,\ldots,A_d)\in Z$ `good' if $A_1x_1,\ldots,A_dx_d$ 
are linearly independent with probability 1.
In \eqref{eq:allchoicesgood}, we obtained that the set of not good elements of 
$Z$ is of measure zero.

Thus, we may choose a sequence $(A_1^{(j)}, A_2^{(j)}, \ldots,A_d^{(j)})$, 
$j=1,2,\ldots$ in $Z$, such that $\|A_i^{(j)}-\id_d\|<1/j$ for all $i$ and $j$, 
and 
$(A_1^{(j)},\ldots,A_d^{(j)})$ is good for each $j$.

Note that for any $j$,
\begin{equation}\label{eq:dominationa}
\left[\det \left(A_1^{(j)}x_1,\ldots, 
A_d^{(j)}x_d\right)\right]^2\leq
|A_1^{(j)}x_1|^2|A_2^{(j)}x_2|^2\ldots 
|A_d^{(j)}x_d|^2,
\end{equation}
and
\begin{equation}\label{eq:dominationb}
\E\left[|A_1^{(j)}x_1|^2|A_2^{(j)}x_2|^2\ldots 
|A_d^{(j)}x_d|^2\right]=d^d.
\end{equation}

\noindent We conclude that

\begin{equation*}
\E\left(\left[\det \left(x_1,\ldots, x_d\right)\right]^2\right)=
\end{equation*}
\begin{equation*}
\E\left(\left[\det \lim_{j\to\infty}\left(A_1^{(j)}x_1,\ldots, 
A_d^{(j)}x_d\right)\right]^2\right)\stackrel{(a)}{=}
\end{equation*}
\begin{equation*}
\E\left(\left[\lim_{j\to\infty}\det \left(A_1^{(j)}x_1,\ldots, 
A_d^{(j)}x_d\right)\right]^2\right)\stackrel{(b)}{=}
\end{equation*}
\begin{equation*}
	\lim_{j\to\infty}\E\left(\left[\det \left(A_1^{(j)}x_1,\ldots, 
	A_d^{(j)}x_d\right)\right]^2\right),
\end{equation*}
where, in (a), we use that the determinant is continuous. In (b),
Lebesgue's Dominated Convergence Theorem may be applied by 
\eqref{eq:dominationa} and \eqref{eq:dominationb}. 

Fix $j$ and let $y_1=A_1^{(j)}x_1,\ldots, y_d=A_d^{(j)}x_d$. 
{In order to emphasize that the assumption 
\eqref{eq:isotropiccenterofmass} 
is not needed, and also for completeness, we repeat Pivovarov's argument.}
For $k=1,\ldots,d-1$, let $P_k$ denote the orthogonal projection of $\Red$ onto 
the linear subspace 
$\spann\{y_1,\ldots,y_k\}^{\perp}$. Thus,
\begin{equation}\label{eq:determinantproduct}
|\det(y_1,\ldots,y_d)|=|y_1||P_1y_2|\cdots|P_{d-1}y_d|.
\end{equation}
Note that with probability $1$, $\rank P_k=d-k$. It follows 
from 
\eqref{eq:isotropictesnor} that $\E|P_ky_{k+1}|^2=d-k$. Fubini's Theorem 
applied 
to \eqref{eq:determinantproduct} completes the proof of 
Lemma~\ref{mainmainlemma}.

\section{Proofs of Theorem~\ref{main} and 
Proposition~\ref{mainProbability}}

Let  $u_1,\ldots u_m\in S^{d-1}$ be a set of vectors satisfying 
\eqref{isotropy} with some positive weights $c_1,\ldots, c_m$. We set the 
probability of each vector $u_i$, $i=1,\ldots m$ as $p_i=c_i/d$, and obtain a 
discrete probability distribution.

Let $u_{i_1},\ldots, u_{i_{d}}$ be independent random vectors from the set 
$u_1,\ldots, u_m$ chosen (with possible repetitions) according to the above 
probability distribution. 

By Lemma~\ref{mainmainlemma}, we have that
\[\E \left ([\det (u_{i_1},\ldots, u_{i_d})]^2\right )=\frac{d!}{d^d}.\]
Since the probability that the random vectors $u_{i_1},\ldots, u_{i_d}$ are 
linearly dependent is positive, 
\[\max [\det (u_{i_1},\ldots, u_{i_d})]^2>\frac{d!}{d^d}.\]
{Our goal is to quantify this inequality by bounding from below the 
probability that the determinant is $0$.
Let 
\[M^2:=\max [\det (u_{i_1},\ldots, u_{i_d})]^2.\]

\noindent Note that if an element of $\{u_1,\ldots, u_m\}$ is selected at least 
twice, 
then $\det (u_{i_1},\ldots,\linebreak[0] u_{i_d})=0$. 
Thus,

\begin{equation*}
\E \left ([\det (u_{i_1},\ldots, u_{i_d})]^2\right )
\leq M^2P_1,
\end{equation*}
} 
where $P_1$ denotes the probability that all indices are pairwise distinct. 
Therefore, 
\[M^2\geq\frac{d!}{d^d}\cdot \frac 1 {P_1}.\]

Note that $P_1$ is a degree $d$ elementary symmetric function of the variables 
$p_1,\ldots, p_m$. Furthermore, $p_1+\ldots+p_m=1$ and $p_i\geq 0$ for all 
$i=1,\ldots, m$.
It can easily be seen (using Lagrange multipliers, {or by induction on $m$}) 
that for fixed $m$ and $d$, 
the maximum of $P_1$ is attained when $p_1=\ldots=p_m=1/m$. Thus,
\[P_1\leq d!\binom{m}{d}\frac{1}{m^d}.\]
In summary, 
\[M^2\geq
\frac{d!}{d^d}\cdot \frac{m^d}{d!}\binom{m}{d}^{-1}
.\]

First, we note that $\gamma(d,m):=\frac{m^d}{{d!}}\binom{m}{d}^{-1}$ is 
decreasing in 
$m$.
Thus, by \eqref{eq:GruberSchuster}, we may assume that $m$ is 
as large as possible, that is, $m=\frac{d(d+1)}{2}$ 
proving the first part of Theorem~\ref{main}.

\subsection{Proof of Theorem~\ref{main} (\ref{main1})}

Let 
{$\gamma(d):=\gamma\left(d,\dbound\right)$.}
We show that $\gamma(d)$ is increasing in $d$. 

With the notation $m:=\dbound$, we note that $(d+1)(d+2)/2=m+d+1$. Thus,

\begin{equation*}
\frac{\gamma(d+1)}{\gamma(d)}=
\frac{(m+d+1)^{d+1}m\cdots (m-d+1)}{m^d (m+d+1)\cdots (m+1)}
=
\frac{(m+d+1)^{d}}{m^d}\cdot
\frac{m\cdots (m-d+1)}{(m+d)\cdots (m+1)}
\end{equation*}

\noindent Thus, we need to show that

\begin{equation*}
1+\frac{d+1}{m}>\sqrt[d]{\left(1+\frac{d}{m}\right)\left(1+\frac{d}{m-1}
\right)\cdots\left(1+\frac{d}{m-d+1}\right)},
\end{equation*}
\noindent which, by the AM/GM inequality follows, if

\begin{equation*}
1+\frac{d+1}{m}\geq 
1+d\frac{\frac{1}{m}+\frac{1}{m-1}+\ldots+\frac{1}{m-d+1}}{d},
\end{equation*}
\noindent which is equivalent to

\begin{equation*}
\frac{d}{m}\geq\frac{1}{m-1}+\frac{1}{m-2}+\ldots+\frac{1}{m-d+1}.
\end{equation*}
For this to hold, it is sufficient to show that for every integer or half of an 
integer $1\leq i\leq d/2$, we have that

\begin{equation}\label{eq:mdinequality}
\frac{2d}{(d-1)m}\geq\frac{1}{m-i}+\frac{1}{m-d+i}.
\end{equation}
After substituting $m=\dbound$, it is easy to see that \eqref{eq:mdinequality} 
holds.

\noindent Finally, $\lim_{d\to\infty}\gamma(d)=e$ follows from Stirling's 
formula.

\subsection{Proof of Theorem~\ref{main} (\ref{csvaltozat}) and (\ref{main3})}
Stirling's formula yields both claims.

\subsection{Proof of Proposition~\ref{mainProbability}}
Let $X$ denote the random variable $X:=[\det (x_1,\ldots, x_d)]^2$, 
$E:=\E(X)=\frac{d!}{d^d}$, and
$q:=\PP\left(X\geq \frac{\lambda E}{P_1}\right)$, where, as in 
the proof of Theorem~\ref{main}, $P_1:=\PP(x_1,\ldots,x_d \text{ are pairwise 
distinct})$.

In the proof of Theorem~\ref{main}, we established
\begin{equation}\label{eq:proofofmain}
P_1\leq (\gamma(d,\mm))^{-1},\mbox{ and thus, } q\leq
\PP\left([\det (x_1,\ldots, x_d)]^2 \geq \lambda\gamma(d,\mm)\cdot 
\frac{d!}{d^d}\right).
\end{equation}

Using the fact that $X$ is at most one, we have
\begin{equation*}
E\leq \frac{\lambda E}{P_1} \PP\left(X< \frac{\lambda E}{P_1} \text{ and } 
x_1,\ldots,x_d \text{ are pairwise distinct}\right)+\PP\left(X\geq 
\frac{\lambda 
E}{P_1}\right).
\end{equation*}

That is, $E\leq \frac{\lambda E}{P_1} (P_1-q)+q$,
and thus, by \eqref{eq:proofofmain}
\begin{equation*}
q\geq \frac{(1-\lambda)E}{1-\frac{\lambda E}{P_1}}\geq
\frac{(1-\lambda)d!}{d^d-{\lambda \gamma(d,\mm)d!}}\geq
(1-\lambda)e^{-d},
\end{equation*}
completing the proof of Proposition~\ref{mainProbability}.

\section{Proof of Theorem~\ref{thm:GS}}\label{sec:GSproof}

First, observe that \eqref{isotropy} holds with some positive scalars $c_i$, 
if and only if, the matrix $\id_d/d$ is in the convex hull of the set 
${\mathcal A}=\{v_i\otimes v_i\st i=1,\ldots,m\}$ in the 
real vector space of $d\times d$ matrices. The set ${\mathcal A}$ is 
contained in the subspace of symmetric matrices with trace $1$, which is of 
dimension 
$d(d+1)/2-1$. Carath\'eodory's theorem \cite{Sch14}*{Theorem~1.1.4} now yields 
the desired upper bound on 
$m$.
	
In the case when both \eqref{isotropy} and \eqref{centeredness} are assumed, we 
lift our vectors into $\R^{d+1}$ as follows. 
Let $\hat v_i=\sqrt{\frac{d}{d+1}}(v_i,1/\sqrt{d})\in\R^{d+1}$. It is easy to 
check that $|\hat v_i|=1$, and that {
\eqref{isotropy} holds for the vectors $\hat v_i$ with some positive 
scalars $\hat c_i$ if, and only if,
\eqref{isotropy} and \eqref{centeredness} 
hold for the vectors $v_i$ with scalars $c_i=\frac{d}{d+1}\hat c_i$. 
}
Now, $\hat{v}_i\otimes\hat{v}_i$, $i=1,\ldots, m$ are 
symmetric $(d+1)\times(d+1)$ matrices of trace one, and their $(d+1,d+1)$th 
entry is $1/(d+1)$. The dimension of this subspace of $\R^{(d+1)\times(d+1)}$ 
is $d(d+3)/2-1$, thus, again, by Carath\'eodory's theorem, the proof is 
complete.

\section{Acknowledgements}
F. Fodor  and T. Zarn\'ocz are supported in part by Hungarian National 
Research, 
Development and Innovation Office NKFIH grant K 116451. 

M. Nasz{\'o}di was partially supported by the
National Research, Development and Innovation Office (NKFIH) grant 
NKFI-K119670 and by the
\'UNKP-17-4 New National Excellence Program of the Ministry of Human 
Capacities.

\end{document}